\documentclass[aap,MSNbibl,citesort,dvips]{arximspdf}

\doi{10.1214/11-AAP768}
\volume{21}
\issue{6}
\pubyear{2011}
\firstpage{2424}
\lastpage{2446}

\makeatletter
\newtheorem{theorem}{Theorem}
\newtheorem{corollary}[theorem]{Corollary}
\newproclaim{definition}[theorem]{Definition}
\newtheorem{lemma}[theorem]{Lemma}
\newtheorem{proposition}[theorem]{Proposition}
\newproclaim{remark}[theorem]{Remark}

\newcommand{\eqref}[1]{(\ref{#1})}
\newcommand{\energy}{\mathcal{E}}
\newcommand{\diver}{\operatorname{div}}
\let\epsilon\varepsilon
\makeatother

\begin{document}
\begin{frontmatter}

\title{Anomalous dissipation in a stochastic inviscid dyadic model}
\runtitle{Anomalous dissipation in a stochastic dyadic model}

\begin{aug}
\author[A]{\fnms{David} \snm{Barbato}\corref{}\ead[label=eb]{barbato@math.unipd.it}},
\author[B]{\fnms{Franco} \snm{Flandoli}\ead[label=ef]{flandoli@dma.unipi.it}}
and
\author[C]{\fnms{Francesco} \snm{Morandin}\ead[label=em]{francesco.morandin@unipr.it}}
\runauthor{D. Barbato, F. Flandoli and F. Morandin}
\affiliation{University of Padova, University of Pisa and University of Parma}
\address[A]{D. Barbato\\
Dipartimento di Matematica Pura\\
\quad e Applicata\\
University of Padova\\
Via Trieste, 63\\
35121 Padova\\
Italy\\
\printead{eb}} 
\address[B]{F. Flandoli\\
Dipartimento di Matematica Applicata\\
University of Pisa\\
Via Buonarroti, 1\\
56127 Pisa\\
Italy\\
\printead{ef}}
\address[C]{F. Morandin\\
Dipartimento di Matematica\\
University of Parma\\
Viale Usberti, 53/A\\
43124 Parma\\
Italy\\
\printead{em}}
\end{aug}

\received{\smonth{12} \syear{2010}}
\revised{\smonth{2} \syear{2011}}

%
\begin{abstract}
A stochastic version of an inviscid dyadic model of turbulence, with
multiplicative noise, is proved to exhibit energy dissipation in spite
of the formal energy conservation. As a consequence, global regular
solutions cannot exist. After some reductions, the main tool is the
escape bahavior at infinity of a certain birth and death
process.\looseness=1
\end{abstract}

%


\vspace*{15pt}
\end{frontmatter}

\section{Introduction}
The dyadic model of turbulence has been introduced in, among others,
\cite{KatPav,Wal,FriPav} and \cite{DesNov}, as a
simplified model of fluid dynamics equations in order to investigate a
number of properties which are out of reach at present for more
realistic models. In this paper we study a suitable random
perturbation of the classical dyadic model under which we are able to
prove anomalous dissipation of energy.

On a complete filtered probability space $ ( \Omega,F_{t},P ) $,
let $ ( W_{n} ) _{n\geq1}$ be a sequence of independent Brownian
motions. Consider the infinite system of stochastic differential
equations in
Stratonovich form
%
\begin{equation}\label{eq_stratonovich}
dX_{n}\,{=}\,( k_{n-1}X_{n-1}^{2}\,{-}\,k_{n}X_{n}X_{n+1} ) \,dt\,{+}\,k_{n-1}%
X_{n-1}\,{\circ}\,\,dW_{n-1}-k_{n}X_{n+1}\,{\circ}\,\,dW_{n}\hspace*{-25pt}
\end{equation}
for $n\geq1$, with $X_{0}(t)=0$. Denote by $l^{2}$ the Hilbert space
of real
square summable sequences $x= ( x_{n} ) _{n\geq1}$ and set
$ \Vert x \Vert^{2}=\sum_{n=1}^{\infty}x_{n}^{2}$. We call
\textit{energy} of $X ( t ) := ( X_{n}(t) ) _{n\geq1}$
the quantity $\mathcal{E} ( t ) :=\frac{1}{2} \Vert X (
t ) \Vert^{2}$. Assume for simplicity to have a deterministic
initial condition
\[
X(0)=x,\qquad  x= ( x_{n} ) _{n\geq1}\in l^{2}.
\]
The sequence of positive numbers $ ( k_{n} ) _{n\geq1}$ will
be specified later on; the most natural case in analogy with fluid
dynamics is $k_{n}=\lambda^{n}$ for some $\lambda>1$.

System (\ref{eq_stratonovich}) is \textit{formally} energy
preserving. By the
Stratonovich form of It\^{o}'s formula (see \cite{Kun}), we have
\begin{eqnarray*}
\tfrac{1}{2}\,dX_{n}^{2} & =&X_{n}\circ \,dX_{n}\\
&=& ( k_{n-1}X_{n-1}^{2}
X_{n}-k_{n}X_{n}^{2}X_{n+1} ) \,dt\\
&&{} +k_{n-1}X_{n-1}X_{n}\circ\,dW_{n-1}-k_{n}X_{n}X_{n+1}\circ\,dW_{n}.
\end{eqnarray*}
If we sum formally these identities and use the boundary condition
$X_{0}(t)=0$, we readily have $\frac{1}{2}d\sum_{n=1}^{\infty}X_{n}^{2}=0$,
namely,
\[
\mathcal{E} ( t ) =\mathcal{E} ( 0 ) ,\qquad
P\mbox{-a.s.}%
\]

The aim of this paper is to prove rigorously an opposite statement, a~property
that we could call \textit{anomalous dissipation}. We need the notion of
energy controlled solution that will be given in the next
section.

\begin{theorem}
\label{main_theorem}Assume $k_{n}=\lambda^{n}$ for some $\lambda>1$. Given
$x\in l^{2}$, let $X ( t ) $ be the unique energy controlled
solution of equation~\eqref{eq_stratonovich}. Then, for all $t>0$,
\[
P\bigl(\mathcal{E}(t)=\mathcal{E}(0)\bigr)<1
\]
and for all $\epsilon>0$ there exist $t$ such that
\[
P\bigl(\mathcal{E}(t)<\epsilon\bigr)>0.
\]
Moreover, if $\mathcal{E} ( 0 ) $ is sufficiently small, then
$\mathcal{E}(t)$ decays to zero at least exponentially fast both
almost surely
and in $\mathcal{L}^{1}$.
\end{theorem}

As a consequence of this theorem, in Section~\ref{sec:lack_reg} we
will also prove that \emph{global regular solutions cannot exist}.

The proof of the theorem is built along Sections 3--6 and concluded
in Section~\ref{proof main theorem}.

Being inspired by fluid dynamics, the results of Theorem~\ref{main_theorem}
could be interpreted as a form of \textit{turbulent dissipation}.
Dynamically speaking, it is a~dissipation due to a very fast cascade
mechanism; energy moves faster and faster from low to high wave
numbers $n$ and escapes to infinity in finite time.

Results of anomalous dissipation for \textit{linear} stochastic
systems with
additive noise have been proved by \cite{MSV1,MSV2}. The
notion of
anomalous dissipation of these papers is different and based on invariant
measures, but conceptually the question is the same. The dyadic model
of the
present paper is, to our knowledge, the first nonlinear stochastic case where
anomalous dissipation is proved to occur. Moreover, the proof is entirely
different from those of the additive noise linear stochastic case.

Nonlinear models with anomalous dissipation have been discovered before
in the
deterministic case (see \cite{Ces,CFP1,CFP2,FriPav,KatPav,KiZlat,Wal,BFM,BFM3}). Our model
is a multiplicative random perturbation of models of these forms.
However, let
us stress that the proof given here is totally different from the
proofs of
the deterministic literature, based on monotonicity and positivity properties
of the deterministic part of system (\ref{eq_stratonovich}). These properties
are lost in the stochastic case.

Theorem \ref{main_theorem} remains true when $k_{n}\leq Cn^{\alpha}$
for some $\alpha>1$. However, we make here the assumption
$k_{n}=\lambda^{n}$ in analogy with the deterministic literature on
dyadic models.

The proof is based on three main ingredients: (i) Girsanov's transform
allows us to reduce the problem to a linear stochastic equation; (ii)
second moments of components satisfy a closed system, without moments
of products; (iii) this closed system is the forward equation of a
birth and death process, the escape of which at infinity can be
understood.

Some issues in this procedure are not trivial. One of them is the
equivalence of laws on infinite time horizon (see
Proposition~\ref{prop:novikov_infty_time}). Its proof is nonstandard;
moreover, it is restricted to a range of values of parameters, the
generalization being open. Concerning the idea that square moments
could satisfy a closed equation, related to jump process, we have been
inspired by previous works (\cite{AirMall,CruzFM,CruzMal} and \cite{MalRen}), however, devoted to different models;
the link here with the nonlinear model and the stochasticity is more
transparent via Girsanov's transform.

\subsection{The multiplicative noise in Euler's equations\label{motivation}}

The noise we introduced in equation~\eqref{eq_stratonovich}, motivated
by energy conservation, may appear peculiar from the physical point of
view. Nevertheless, it is the natural choice if we compare the dyadic
model with some other equations of fluid dynamics like Euler's equations
or diffusions of passive scalars.

Let us give a picture of this analogy in the case of Euler's equations,
which have the form
\[
\frac{\partial u}{\partial t}+u\cdot\nabla u+\nabla p=0,\qquad  \diver u=0
\]
with appropriate initial and boundary conditions ($u$ and $p$ are the velocity
and pressure field, resp.). Let us think of periodic boundary
conditions for sake of simplicity. The Lagrangian motion of particles
is given
by the equation
\[
\frac{dY ( t ) }{dt}=u ( t,Y ( t ) ) .
\]
A natural way to randomly perturb Euler dynamics (see \cite{MR1}) is
by adding a white noise to the Lagrangian motion,
\[
dY ( t ) =u ( t,Y ( t ) ) \,dt+\sum_{j}%
\sigma_{j} ( t,Y ( t ) ) \,dW^{j} ( t ),
\]
where $ ( W^{j} ( t ) ) _{t\geq0}$ are independent
Brownian motions and $\sigma_{j}$ are given vector fields. By standard rules
of stochastic calculus, applied formally, one can see that Euler equations
take the stochastic form (with a new pressure~$\widetilde{p}$)
%
\begin{equation}\label{stoch euler eq}
du+ [ u\cdot\nabla u+\nabla\widetilde{p} ] \,dt+\sum_{j}\sigma
_{j}\cdot\nabla u\circ \,dW^{j} ( t ) =0, \qquad \diver u=0,
\end{equation}
where Stratonovich operation has to be used. Rigorous results and
physical arguments in support of this kind of stochastic perturbation
of the Lagrangian motion and the corresponding PDE with multiplicative
Stratonovich noise (in the viscous case) can be found in \cite{MR1,MR2}. In addition, let us mention the wide literature on
stochastic passive scalar equations (see, e.g., \cite{LR}) where
multiplicative Stratonovich noise of the form above is used.

In abstract form, equation (\ref{stoch euler eq}) takes the form
%
\begin{equation}\label{abstract stoch Euler eq}
du+B ( u,u ) \,dt+B ( \circ \,dW,u ),
\end{equation}
where $W ( t ) :=\sum_{j}\sigma_{j} ( x ) W^{j} (
t ) $ and $B(u,v):=u\cdot\nabla v$. Notice that for sufficient regular
$u$ and $v$ the following identity holds $ \langle B(u,v),v \rangle
=0$.

In \cite{CFP1} the authors argued that after Fourier or wavelets transforms and
certain simplifications, the deterministic system
%
\begin{equation}\label{dyadic_det}%
\frac{dX_{n}}{dt}=k_{n-1}X_{n-1}^{2}-k_{n}X_{n}X_{n+1},\qquad  n\geq1,
\end{equation}
describes some idealized features of the deterministic equation $\frac
{du}%
{dt}+B ( u,u ) =0$. Equation (\ref{dyadic_det}) has the form
$\frac{d}{dt}X=\tilde{B}(X,X)$ where
\[
\tilde{B}(X,Y)_{(n)}=k_{n-1}X_{n-1}Y_{n-1}-k_{n}X_{n}Y_{n+1}.
\]
Formally we have $ \langle\tilde{B}(X,Y),Y \rangle=0$. Thanks to
$ \langle B(u,u),u \rangle=0$ and $ \langle\tilde
{B}(X,X),\allowbreak X \rangle=0$ the perturbation of $u$ and $X$ does not modify the
energy balance (at a formal level).

Standing this idealized discretization $\tilde{B}(X,Y)$ of $B(u,v)$, the
natural analog of equation (\ref{abstract stoch Euler eq}) is
\[
dX_{n}+\tilde{B}(X,X)_{(n)}\,dt+\tilde{B}(\circ\, dW,X)_{(n)},\qquad
n\geq1,
\]
which is precisely system (\ref{eq_stratonovich}).

\section{It\^{o}'s formulation}

For the rigorous formulation of equation (\ref{eq_stratonovich}) and a~basic
theorem of existence and uniqueness, we follow \cite{BFM2}. The It\^
{o} form
of equation (\ref{eq_stratonovich}) is
\begin{eqnarray}\label{eq_Ito}
dX_{n} & =& ( k_{n-1}X_{n-1}^{2}-k_{n}X_{n}X_{n+1} ) \,dt+k_{n-1}%
X_{n-1}\,dW_{n-1}\nonumber\\[-8pt]\\[-8pt]
&&{} -k_{n}X_{n+1}\,dW_{n}-\tfrac{1}{2} ( k_{n}^{2}+k_{n-1}^{2} ) X_{n}\,dt.\nonumber
\end{eqnarray}
Let us define the concept of weak solution for this equation. By a filtered
probability space $ ( \Omega,F_{t},P ) $ we mean a probability
space $ ( \Omega,F_{\infty},P ) $ and a~right-continuous filtration
$ ( F_{t} ) _{t\geq0}$ such that $F_{\infty}$ is the $\sigma
$-algebra generated by $\bigcup_{t\geq0}F_{t}$.\vspace*{-2pt}

\begin{definition}
\label{def sol 1}Given $x\in l^{2}$, a weak solution of equation
(\ref{eq_stratonovich}) in $l^{2}$ is a filtered probability space $ (
\Omega,F_{t},P ) $, a sequence of independent Brownian motions $ (
W_{n} ) _{n\geq1}$ on $ ( \Omega,F_{t},P ) $ and an $l^{2}%
$-valued stochastic process $ ( X_{n} ) _{n\geq1}$ on $ (
\Omega,F_{t},P ) $ with continuous adapted components $X_{n}$, such
that
\begin{eqnarray*}
X_{n} ( t ) & =&x_{n}+\int_{0}^{t} \bigl( k_{n-1}X_{n-1}%
^{2} ( s ) -k_{n}X_{n} ( s ) X_{n+1} ( s )
\bigr) \, ds\\
&&{} +\int_{0}^{t}k_{n-1}X_{n-1} ( s ) \,dW_{n-1} ( s )
-\int_{0}^{t}k_{n}X_{n+1} ( s ) \,dW_{n} ( s ) \\
&&{} -\int_{0}^{t}\frac{1}{2} ( k_{n}^{2}+k_{n-1}^{2} ) X_{n} (
s ) \,ds
\end{eqnarray*}
for each $n\geq1$, with $X_{0}=0$. We denote this solution by
\[
( \Omega,F_{t},P,W,X )
\]
or simply by $X$.\vspace*{-2pt}
\end{definition}

\begin{definition}
\label{def sol 2}We call energy controlled solutions the solutions of
Definition~\ref{def sol 1} which satisfy
%
\begin{equation}\label{energy_inequality}
P \Biggl( \sum_{n=1}^{\infty}X_{n}^{2} ( t ) \leq\sum_{n=1}^{\infty
}x_{n}^{2} \Biggr) =1%
\end{equation}
for all $t\geq0$.
\end{definition}

The following simple proposition (proved in \cite{BFM2}) clarifies
that a
process satisfying (\ref{eq_Ito}) rigorously satisfies also
(\ref{eq_stratonovich}).\vspace*{-2pt}

\begin{proposition}
If $X$ is a weak solution, for every $n\geq1$ the process $ (
X_{n} ( t ) ) _{t\geq0}$ is a continuous semimartingale,
hence, the two Stratonovich integrals
\[
\int_{0}^{t}k_{n-1}X_{n-1} ( s ) \circ \,dW_{n-1} ( s )
-\int_{0}^{t}k_{n}X_{n+1} ( s ) \circ \,dW_{n} ( s )
\]
are well defined and equal to
\begin{eqnarray*}
&& \int_{0}^{t}k_{n-1}X_{n-1} ( s ) \,dW_{n-1} ( s )
-\frac{1}{2}\int_{0}^{t}k_{n-1}^{2}X_{n} ( s ) \,ds\\
&&\qquad{} -\int_{0}^{t}k_{n}X_{n+1} ( s ) \,dW_{n} ( s )
-\frac{1}{2}\int_{0}^{t}k_{n}^{2}X_{n} ( s ) \,ds.
\end{eqnarray*}
Hence, $X$ satisfies the Stratonovich equations\vadjust{\goodbreak}
(\ref{eq_stratonovich}).
\end{proposition}

The main result proved in \cite{BFM2} is the well posedness in the
weak probabilistic sense in the class of energy controlled
solutions.\vspace*{-2pt}

\begin{theorem}
Given $ ( x_{n} ) \in l^{2}$, there exists one and only one energy
controlled solution of equation (\ref{eq_stratonovich}).\vspace*{-2pt}
\end{theorem}

\section{Girsanov's transformation}

Formally, let us write equation (\ref{eq_Ito}) in the form\vspace*{-2pt}
\begin{eqnarray*}
dX_{n} & =&k_{n-1}X_{n-1} ( X_{n-1}\,dt+dW_{n-1} ) -k_{n}%
X_{n+1} ( X_{n}\,dt+dW_{n} ) \\[-2pt]
&&{} -\tfrac{1}{2} ( k_{n}^{2}+k_{n-1}^{2} ) X_{n}\,dt.\vspace*{-2pt}
\end{eqnarray*}
The simple idea is that $X_{n}\,dt+dW_{n}$ is a Brownian motion with
respect to
a~new measure $Q$ on $ ( \Omega,F ) $, simultaneously for every
$n$, hence, the equations become linear SDEs under $Q$. We use details about
Girsanov's theorem that can be found in \cite{RevuzYor}, Chapter VIII,
and an
infinite dimensional version proved in \cite{Bensou,Kozlov,DZ}.

Assume that $ ( X_{n} ) _{n\geq1}$ is an energy controlled
solution. Due to the boundedness of $\sum_{n=1}^{\infty}X_{n}^{2} (
t ) $ [see (\ref{energy_inequality})], the process $Y_{t}:=-\sum
_{n=1}^{\infty}\int_{0}^{t}X_{n} ( s ) \,dW_{n} ( s ) $ is
well defined, is a martingale and its quadratic variation $ [ Y,Y ]
_{t}$ is\break $\int_{0}^{t}\sum_{n=1}^{\infty}X_{n}^{2} ( s )\, ds$. For
the same reason, Novikov criterium applies, so\break $\mathcal{N} ( Y )
_{t}:=\exp( Y_{t}- [ Y,Y ] _{t} ) $ is a strictly
positive martingale. Define the set function $Q$ on $\bigcup_{t\geq
0}F_{t}$ by
setting\vspace*{-2pt}
%
\begin{equation} \label{QT}\qquad
\frac{dQ}{dP} \bigg\vert_{F_{t}}=\mathcal{N} ( Y )
_{t}=\exp\Biggl( -\sum_{n=1}^{\infty}\int_{0}^{t}X_{n} ( s )
\,dW_{n} ( s ) -\frac{1}{2}\int_{0}^{t}\sum_{n=1}^{\infty}X_{n}%
^{2} ( s ) \,ds \Biggr)\vspace*{-2pt}
\end{equation}
for every $t\geq0$. We also denote by $Q$ its extension to the terminal
$\sigma$-field~$F_{\infty}$. In general we cannot prove it is absolutely
continuous with respect to~$P$, but we shall see at least a case when
this is
true. Notice also that $Q$ and~$P$ are equivalent on each $F_{t}$, by the
strict positivity. Define\vspace*{-2pt}
\[
B_{n} ( t ) =W_{n} ( t ) +\int_{0}^{t}X_{n} (
s ) \,ds.\vspace*{-2pt}
\]
Under $Q$, $ ( B_{n} ( t ) ) _{n\geq1,t\in[
0,T ] }$ is a sequence of independent Brownian motions. Since\vspace*{-2pt}
\begin{eqnarray*}
\int_{0}^{t}k_{n-1}X_{n-1} ( s ) \,dB_{n-1} ( s ) &
=&\int_{0}^{t}k_{n-1}X_{n-1} ( s ) \,dW_{n-1} ( s ) \\[-2pt]
&&{} +\int_{0}^{t}k_{n-1}X_{n-1} ( s ) X_{n-1} ( s ) \,ds\vspace*{-2pt}
\end{eqnarray*}
and similarly for $\int_{0}^{t}k_{n}X_{n+1} ( s ) \, dB_{n} (
s ) $, we see that\vspace*{-2pt}
\begin{eqnarray*}
X_{n} ( t ) & =&X_{n} ( 0 ) +\int_{0}^{t}k_{n-1}%
X_{n-1} ( s ) \,dB_{n-1} ( s ) -\int_{0}^{t}k_{n}%
X_{n+1} ( s ) \,dB_{n} ( s ) \\[-2pt]
&&{} -\int_{0}^{t}\frac{1}{2} ( k_{n}^{2}+k_{n-1}^{2} ) X_{n} (
s )\, ds.\vadjust{\goodbreak}
\end{eqnarray*}
This is a linear stochastic equation. Girsanov's transformation has
removed the
nonlinearity. Let us collect the previous facts.

\begin{theorem}
\label{teo_Girsanov}If $ ( \Omega,F_{t},P,W,X ) $ is an energy
controlled solution of the nonlinear equation (\ref
{eq_stratonovich}), then it
satisfies the linear equation
%
\begin{equation}
\label{eq:stoch_linear_Q}dX_{n}=k_{n-1}X_{n-1}\,dB_{n-1}-k_{n}X_{n+1}%
\,dB_{n}-\tfrac{1}{2} ( k_{n}^{2}+k_{n-1}^{2} ) X_{n}\,dt,
\end{equation}
where the processes
\[
B_{n} ( t ) =W_{n} ( t ) +\int_{0}^{t}X_{n} (
s ) \,ds
\]
are a sequence of independent Brownian motions on $ ( \Omega
,F_{t},Q ) $, $Q$ defined by~(\ref{QT}).
\end{theorem}

One may also check that
\[
dX_{n}=k_{n-1}X_{n-1}\circ \,dB_{n-1}-k_{n}X_{n+1}\circ \,dB_{n}%
\]
so the previous computations could be described at the level of
Stratonovich calculus.

\section{Closed equation for $\mathbb{E}^{Q}[X_{n}^{2}(t)]$}

Let $ ( \Omega,F_{t},P,W,X ) $ be an energy controlled solution of
the nonlinear equation (\ref{eq_stratonovich}) with initial condition
$x\in
l^{2}$ and let $Q$ be the measure given by Theorem \ref{teo_Girsanov}. Denote
by $\mathbb{E} ^{Q}$ the mathematical expectation on $ ( \Omega
,F_{t},Q ) $. We have
\begin{eqnarray*}
\tfrac{1}{2}\,dX_{n}^{2} & =&X_{n}\,dX_{n}+\tfrac{1}{2}\,d [ X_{n} ]
_{t}\\
& =&-\tfrac{1}{2} ( k_{n}^{2}+k_{n-1}^{2} ) X_{n}^{2}\,dt+dM_{n}%
+\tfrac{1}{2} ( k_{n-1}^{2}X_{n-1}^{2}+k_{n}^{2}X_{n+1}^{2} ) \,dt,
\\
\tfrac{1}{2}\,dX_{n}^{4} & =&4X_{n}^{3}\,dX_{n}+\tfrac{12}{2}X_{n}^{2}\,d [
X_{n} ] _{t}\\
& =&-\tfrac{1}{2} ( k_{n}^{2}+k_{n-1}^{2} ) X_{n}^{4}\,dt+dM_{n}%
+\tfrac{12}{2}X_{n}^{2} ( k_{n-1}^{2}X_{n-1}^{2}+k_{n}^{2}X_{n+1}%
^{2} ) \, dt,
\\
de^{X_{n}^{2}}&=&e^{X_{n}^{2}}\,dX_{n}^{2}+\tfrac{1}{2}e^{X_{n}^{2}}\,d [
X_{n}^{2},X_{n}^{2} ],
\\
de^{X_{n}} & =&e^{X_{n}}\,dX_{n}+\tfrac{1}{2}e^{X_{n}}\,d [ X_{n}%
,X_{n} ] \\
& =&\cdots e^{X_{n}} ( k_{n}^{2}+k_{n-1}^{2} ) X_{n}\,dt+dM_{n}+\tfrac
{1}{2}e^{X_{n}}\tfrac{1}{2} ( k_{n-1}^{2}X_{n-1}^{2}+k_{n}^{2}X_{n+1}%
^{2} ) \,dt,
\end{eqnarray*}
where
\[
M_{n} ( t ) =\int_{0}^{t}k_{n-1}X_{n-1} ( s )
X_{n} ( s ) \,dB_{n-1} ( s ) -\int_{0}^{t}k_{n}%
X_{n} ( s ) X_{n+1} ( s ) \, dB_{n} ( s ) .
\]
Notice that
%
\begin{equation}\label{forth_order}
\mathbb{E} ^{Q}\int_{0}^{T}X_{n}^{4} ( t ) \,dt<\infty
\end{equation}
for each $n\geq1$. Indeed, for an energy controlled solution, from
(\ref{energy_inequality}) we have, with $P$-probability one,
\[
\sum_{n=1}^{\infty}X_{n}^{4} ( t ) \leq\max_{n}X_{n}^{2} (
t ) \sum_{n=1}^{\infty}X_{n}^{2} ( t ) \leq\Biggl( \sum
_{n=1}^{\infty}x_{n}^{2} \Biggr) ^{2}.
\]
But $P$ and $Q$ are equivalent on $F_{t}$, hence,
\[
Q \Biggl( \sum_{n=1}^{\infty}X_{n}^{4} ( t ) \leq\Biggl( \sum
_{n=1}^{\infty}x_{n}^{2} \Biggr) ^{2} \Biggr) =1.
\]
This implies (\ref{forth_order}).

From (\ref{forth_order}), $M_{n} ( t ) $ is a martingale for each
$n\geq1$. Moreover,\break $\mathbb{E} ^{Q} [ \sum_{n=1}^{\infty}X_{n}%
^{2} ( t ) ] <\infty$ because $ ( X_{n} )
_{n\geq1}$ is an energy controlled solution [again, as above, condition
(\ref{energy_inequality}) is invariant under the change of measure
$P\leftrightarrow Q$ on~$F_{t}$] and thus, in particular, $\mathbb{E}
^{Q} [ X_{n}^{2} ( t ) ] $ is finite for each $n\geq1$.
From the previous equation we deduce the following.

\begin{proposition}
\label{prop valori medi}For every energy controlled solution $X$,
$\mathbb{E}
^{Q} [ X_{n}^{2} ( t ) ] $ is finite for each $n\geq1$
and satisfies
\begin{eqnarray*}
\frac{d}{dt}\mathbb{E} ^{Q} [ X_{n}^{2} ] & =&- ( k_{n}%
^{2}+k_{n-1}^{2} ) \mathbb{E} ^{Q} [ X_{n}^{2} ] \\
&&{} +k_{n-1}^{2}\mathbb{E} ^{Q} [ X_{n-1}^{2} ] +k_{n}^{2}\mathbb{E}
^{Q} [ X_{n+1}^{2} ]
\end{eqnarray*}
for $t\geq0$.
\end{proposition}

The first remarkable fact of this result is that $\mathbb{E} ^{Q} [ X_{n}
^{2} ] $ satisfies a closed equation. The second one is that this is the
forward equation of a continuous-time Markov chain, as we shall discuss
in the
next section. See \cite{AirMall,CruzFM} for different
examples with
the same structure.

\section{Associated birth and death process}

In this section we will make thorough use of birth and death
processes. We do not suppose that all the readers are familiar with the
field, so we will be more detailed.

Let us set
\[
p_{n} ( t ) =\frac{1}{ \| x \| ^{2}}\mathbb{E}
^{Q} [ X_{n}^{2} ( t ) ] ,\qquad  p ( t )
= ( p_{n} ( t ) ) _{n\geq1},\qquad  t\geq0,
\]
and set also $p_{0} ( t ) \equiv0$. Introduce the positive numbers
$(\lambda_{n})_{n\geq1}$ and $(\mu_{n})_{n\geq1}$, defined as
\[
\lambda_{n} =k_{n}^{2}, \qquad \mu_{n} =k_{n-1}^{2}.
\]
By Proposition~\ref{prop valori medi}, we have
%
\begin{equation}
\label{forw_eq}
\quad \cases{
\displaystyle\frac{d}{dt}p_{n} ( t ) =-(\lambda_{n}+\mu_{n})p_{n} (
t ) +\lambda_{n-1}p_{n-1} ( t ) +\mu_{n+1}p_{n+1} (
t ) ,&\quad$ t\geq0$,\cr
\displaystyle p_{n} ( 0 ) =\frac{x_{n}^{2}}{ \| x \| ^{2}}.%
}\hspace*{-10pt}
\end{equation}
We observe that $\sum_{n=1}^{\infty}p_{n} ( t ) =1$ when $t=0$ and,
moreover,
%
\begin{equation}
\label{eq:condit_fwd}\sum_{n=1}^{\infty}p_{n} ( t ) \leq1
\end{equation}
for all $t>0$ (since $X$ is an energy controlled solution).

The system~\eqref{forw_eq} can be conveniently put in matrix form,
$\frac
d{dt}p(t)=p(t)A$, where $A$ is an infinite matrix with null row sums and
nonnegative off-diagonal entries.

In the theory of continuous-time Markov chains this is usually referred
to as
a $q$-matrix. Since it has tridiagonal form, all the processes with
$q$-matrix~$A$ will be \textit{birth and death processes}.

By studying $A$ we will be able to identify exactly one process $\xi
_{t}$ on
some new probability space $(S,\mathcal{S}, \mathcal{P} )$ such that
$p_{n}(t)=\mathcal{P} (\xi_{t}=n)$. Since $\xi$ will turn out to be
\emph{dishonest} (meaning that $\mathcal{P} $-a.s.\ $\xi$ will
escape to
infinity in finite time), the conclusion will be that $\lim
_{t\rightarrow
\infty}\sum_{n=1}^{\infty}\mathbb{E} ^{Q}[X_{n}^{2}(t)]=0$.

\subsection{Minimal process}

In general, given a $q$-matrix $A$, there can be many processes $\chi$ with
different laws $y_{n}(t)=\mathcal{P} (\chi(t)=n)$, all satisfying
either the
forward $y^{\prime}=yA$ or the backward $y^{\prime}=Ay$ equations associated
with~$A$. Whether the solutions of the two systems are unique depends
on some
well-studied properties of the $q$-matrix.

In the present case, $A$ is stable (no $-\infty$ entries appear in the
diagonal) and conservative (no mass disappears at zero because $\mu_{1}=0$).
It is well known that to any stable $q$-matrix is associated a process, called
\emph{minimal}, whose law satisfies both systems of equations.

The latter is the naive process that anyone would construct from $A$, as
follows. Given a probability space $(S,\mathcal{S}, \mathcal{P} )$, let
$\xi_{t}$ be a continuous-time Markov chain on the positive integers, with
initial distribution
\[
\mathcal{P} (\xi_{0}=n)=p_{n}(0),\qquad  n=1,2,\dots,
\]
and jump rates given by $A$ entries, that is, the process waits in a
state $n$
for an exponential time with rate $\lambda_{n}+\mu_{n}$ and then
jumps at
$n+1$ or $n-1$ with probabilities $\pi_{n}$ and $1-\pi_{n}$, respectively,
where
\[
\pi_{n}:=\frac{\lambda_{n}}{\lambda_{n}+\mu_{n}}.
\]
Let $\tau\in[0,\infty]$ denote the first time such that in $[0,\tau
)$ the
process has undergone infinitely many jumps. We say that the process reaches
the boundary at time $\tau$ and we give no special ``return'' rule if the
process reaches the boundary in finite time. Hence, if for $\omega\in S$,
$\tau(\omega)<\infty$, then $\xi_{t}(\omega)$ is not defined for
$t\geq
\tau(\omega)$. Notice that, given $s>0$, $\mathcal{P} (\tau>s)=\sum
_{n=1}^{\infty}\mathcal{P} (\xi_{s}=n)$ could be less than 1.

If the minimal solution of a $q$-matrix is honest, it is the unique solution
for each one of the two systems and the $q$-matrix itself is called regular.
As anticipated, the minimal solution, which is the law of the process
described above, is not regular if the coefficients $k_{n}$ grow too fast
(Proposition~\ref{prop:A_not_regular} below), nevertheless it is the unique
solution of the forward equations (Proposition~\ref{prop:uniq_fwd} below),
while the backward equations have infinite solutions.

This uniqueness is very important because it ensures that $\mathcal{P}
(\xi_{t}=n)=p_{n}(t):=\|x\|^{-2}\mathbb{E} ^{Q}[X_{n}^{2}(t)]$. If we denote
by $\mathcal{E} $ the total energy of $X$,
%
\begin{equation}
\mathcal{E} (t) :=\frac12\sum_{n=1}^{\infty}X_{n}^{2}(t)
\end{equation}
this means in particular that we can study $\mathbb{E} ^{Q}[\mathcal
{E} (t)]$
through $\mathcal{P} (\tau>t)$,
%
\begin{equation}
\label{eq:Qavg_ener_via_tau}\mathbb{E} ^{Q}[\mathcal{E} (t)] =\frac
12\sum_{n=1}^{\infty}\mathbb{E} ^{Q}[X_{n}^{2}(t)] =\mathcal{E}
(0)\sum
_{n=1}^{\infty}p_{n}(t) =\mathcal{E} (0)\mathcal{P} (\tau>t),
\end{equation}
which will be the aim of Section~\ref{sec:tau}.

\begin{proposition}
\label{prop:A_not_regular}
$\!\!\!$The $q$-matrix $A$ is not regular if and
only if
\mbox{$\sum_{n}nk_{n}^{-2}\,{<}\,\infty$}.
\end{proposition}

For the proof we make use of results by Reuter and Anderson, which are
efficiently exposed in the book by the latter \cite{Anderson}. It is, however,
not too difficult an exercise to prove the ``if'' direction with elementary
notions. Truly, Proposition~\ref{prop_tau_exp_bd} and
Lemma~\ref{lemma_less_expl_from_1} below provide such an argument and
we refer
the reader who wants some insight to them.

\begin{pf*}{Proof of Proposition~\ref{prop:A_not_regular}}
By Corollary~2.2.5
of~\cite{Anderson}, in the conservative case, the minimal solution is honest
if and only if the backward equations have a unique solution.

By Theorem~3.2.2 of~\cite{Anderson} the $q$-matrix of a birth and death
process has a unique solution of the backward equations if and only if the
following quantity is infinite:
\[
R=\sum_{n=1}^{\infty} \biggl( \frac1{\lambda_{n}} +\frac{\mu
_{n}}{\lambda_{n}\lambda_{n-1}} +\frac{\mu_{n}\mu_{n-1}}{\lambda_{n}\lambda
_{n-1}\lambda_{n-2}} +\cdots+\frac{\mu_{n}\cdots\mu_{2}}{\lambda
_{n}\cdots\lambda_{2}\lambda_{1}} \biggr) .
\]
Since $\lambda_{n}=\mu_{n+1}$, we get $R=\sum_{n}n\lambda
_{n}^{-1}=\sum
_{n}nk_{n}^{-2}$.
\end{pf*}

\begin{proposition}
\label{prop:uniq_fwd}
The forward system of equations~\eqref{forw_eq},
together with condition~\eqref{eq:condit_fwd} admits a unique solution.
\end{proposition}

Here again the proposition can be seen as a simple application of a result
from the book by Anderson, specifically Theorem~3.2.3 of \cite{Anderson}.

Uniqueness could also be proved with an analytic approach, based on the
parabolic structure of the equation which is apparent if we remember
$\lambda_{n}=\mu_{n+1}$ and we rewrite~\eqref{forw_eq} as
%
\begin{equation}
\label{backw_eq}
\quad \frac{d}{dt}p_{n}(t) =\lambda_{n}\bigl(p_{n+1}(t)-p_{n}%
(t)\bigr)-\lambda_{n-1}\bigl(p_{n}(t)-p_{n-1}(t)\bigr) ,\qquad  t\geq0.
\end{equation}
Nevertheless, we believe it would be interesting to show a completely different
and entirely elementary proof that maybe could also be used when the parabolic
nature is lost and the associated process is no more a simple birth and death.

\begin{pf*}{Proof of Proposition \ref{prop:uniq_fwd}}
By linearity we can suppose $p_{n}(0)=0$, with condition~\eqref{eq:condit_fwd}
still holding and $-1\leq p_{n}(t)\leq1$.

Suppose by contradiction that $(p_{n})_{n\in\mathbb{N}}$ is a nonzero
solution. Without loss of generality we can suppose
$p_{1}(t_{0})=\delta>0$
for some $t_{0}>0$ [take the largest~$n_{0}$ such that $p_{n}\equiv0$
for all
$n<n_{0}$, so we have $p_{n_{0}}(t_{0})>0$ for some $t_{0}>0$; then
shift and
rename the indexes of the sequence $(p_{n})$ in such a way that the new
$p_{1}$ is the old $p_{n_{0}}$]. Define the partial sums
\[
\phi_{n}(t):=\sum_{j=1}^{n}p_{j}(t).
\]
We notice that a simple computation starting from~\eqref{backw_eq} yields
%
\begin{equation}
\label{cumul_eq}\frac{d}{dt}\phi_{n}=\lambda_{n}(p_{n+1}-p_{n}).
\end{equation}
Then define the times
%
\begin{equation}
\label{def_t_n}t_{n}:=\inf\{t|\phi_{n}(t)\geq n\delta\},\qquad  n\geq1.
\end{equation}
We claim that for all $n\geq1$,
\[
t_{n}\leq t_{n-1} \quad \mbox{and}\quad  p_{n}(t_{n})\geq\delta
\]
so that the sequence $(t_{n})_{n\geq0}$ is finite (in fact,
decreasing), in
contradiction with the position $\phi_{n}(t)\leq1$ for all $n$ and
for all $t$.

We shall prove the claim by induction.

For $n=1$, by definition $t_{1}\leq t_{0}<+\infty$ and
$p_{1}(t_{1})=\phi
_{1}(t_{1})\geq\delta$.

Let us suppose that the claim holds for $n$. By the definition of
\mbox{$t_{n},\frac{d}{dt}\phi_{n}(t_{n})\,{\geq}\,0$}.

By~\eqref{cumul_eq}, this implies $p_{n+1}(t_{n})\geq p_{n}(t_{n})\geq
\delta$
and hence
\[
\phi_{n+1}(t_{n})= \phi_{n}(t_{n})+p_{n+1}(t_{n})\geq n\delta+\delta
\]
thus, $t_{n+1}\leq t_{n}$. This implies that $\phi_{n}(t_{n+1})\leq
n\delta$,
so that necessarily
\[
p_{n+1}(t_{n+1})\geq\delta.
\]
The induction and the proof are complete.
\end{pf*}

We remark the fact that given the condition $\sum_{n=1}^{\infty}nk_{n}
^{-2}<\infty$, the forward equations have a unique solution while the
backward have infinitely many. This fact might appear a bit
disconcerting if one notices that $A$ is symmetric and hence, forward
and backward equations are formally identical. The explanation is that
any proper solution $p_n(t)$ of the forward system of equations must
be summable in the sense that $\sum_np_n(t)\leq1$. On the contrary,
if $\{q_n(t;k)\}_n$ is a solution of the backward equations with
initial condition $q_n(0;k)=\delta_{k,n}$, it must satisfy
$\sum_kq_n(t;k)\leq1$ for all $n$.

\subsection{Time of escape}

\label{sec:tau}

In this section we study the law of $\tau$, the time of escape to
infinity of
the minimal process. The main result is Proposition~\ref{prop_tau_exp_bd},
which is generalized by Lemma~\ref{lemma_less_expl_from_1}.

\begin{lemma}
\label{lemma_chain_geom_visit} Suppose $\sum_{i=1}^{\infty
}k_{i}^{-2}<\infty$
and that the minimal process starts from~$1$. For $n\geq1$, the number
of times
the minimal process visits state $n$ is a geometric r.v.\ with mean
$(k_{n}^{2}+k_{n-1}^{2})\sum_{i=n}^{\infty}k_{i}^{-2}$.
\end{lemma}

\begin{pf}
We follow ideas from Feller \cite{Feller}. Let $p_{i,j}$ denote the transition
probabilities of the discrete time Markov chain embedded in continuous-time
minimal process and let $\sigma^{(i)} =\{\sigma_{n}^{(i)}\}_{n>i}$
denote the
probabilities that the chain starting from states $n$ larger than $i$ will
never get to $i$. Then $\sigma^{(i)}$ is the maximal solution of
%
\begin{equation}\label{eq:syst_sigma_k}
x_{n}=\sum_{j>i}p_{n,j}x_{j},\qquad  n>i,
\end{equation}
satisfying $0\leq x_{n}\leq1$ for all $n$. (This solution can be zero.)

If we let $x_{i}=0$ for sake of notation, in our case the system
(\ref{eq:syst_sigma_k}) reduces to
\[
x_{n} =\frac{\mu_{n}}{\lambda_{n}+\mu_{n}}x_{n-1}+\frac{\lambda
_{n}%
}{\lambda_{n}+\mu_{n}}x_{n+1},\qquad n\geq i+1,
\]
yielding
\[
x_{n+1}-x_{n} =\frac{\mu_{n}}{\lambda_{n}}(x_{n}-x_{n-1}),\qquad
n\geq i+1,
\]
and then by induction, for $n\geq i$,
\begin{eqnarray*}
x_{n+1}-x_{n} & =&x_{i+1}\prod_{j=k+1}^{n}\frac{\mu_{j}}{\lambda
_{j}}=x_{i+1}\frac{k_{i}^{2}}{k_{n}^{2}},\\
x_{n} & =&x_{i+1}k_{i}^{2}\sum_{m=i}^{n-1}k_{m}^{-2}.%
\end{eqnarray*}
By hypothesis the sums are bounded and hence, the maximal solution is obtained
by choosing $x_{i+1}$ such that $\lim_{n}x_{n}=1$, that is,
\[
\sigma_{i+1}^{(i)}= \Biggl( k_{i}^{2}\sum_{m=i}^{\infty} k_{m}^{-2}%
\Biggr) ^{-1},
\]
hence, the chain is transient.

Now suppose that the chain is starting from 1. It will visit $i$ at least
once. When it does, the probability that it is the last visit is
$p_{i,i+1}\sigma_{i+1}^{(i)}$, so by strong Markov property, the total number
of visits to $i$ is a geometric random variable with mean
\[
\bigl(p_{i,i+1}\sigma_{i+1}^{(i)} \bigr)^{-1} =(k_{i}^{2}+k_{i-1}^{2}%
)\sum_{m=i}^{\infty}k_{m}^{-2}
\]
as required.
\end{pf}

\begin{proposition}
\label{prop_tau_exp_bd} Suppose $\nu_{\infty}:=\sum_{n=1}^{\infty
}nk_{n}%
^{-2}<\infty$ and that the minimal process starts from $1$. Let
$T_{n}$ be the
total time the minimal process spends in the state $n$,
\[
T_{n}:=\mathcal{L}\{t\geq0\dvtx\xi_{t}= n\},
\]
so that the time of escape at infinity is $\tau=\sum_{n=0}^{\infty}T_{n}$.

Then for all $n\geq1$, $T_{n}$ is an exponential r.v.\ with mean $\nu
_{n}:=\sum_{i=n}^{\infty}k_{i}^{-2}$ and in particular
\[
\mathbb{E} ^{\mathcal{P} }(\tau)=\sum_{n=1}^{\infty}\nu_{n}=\sum
_{n=1}%
^{\infty}nk_{n}^{-2}=\nu_{\infty}.
\]
Moreover, there exists $h>0$ such that for all $t$
\[
e^{-t/\nu_{1}} \leq\mathcal{P} (\tau>t) \leq e^{-t/\nu_{\infty}+h}.
\]

\end{proposition}

\begin{pf}
The total time spent in a state $n$ is the sum of many i.i.d.\ exponential
waiting times of rate $k_{n}^{2}+k_{n-1}^{2}$.

Since the sum of a geometric number of i.i.d.\ exponential r.v.'s is
exponential, by Lemma~\ref{lemma_chain_geom_visit}, $T_{n}$ is
exponential and
its mean is as required.

The lower bound comes easily from $\tau=\sum_{n}T_{n}$, since
\[
\mathcal{P} (\tau>t) \geq\mathcal{P} (T_{1}>t) =e^{-t/\nu_{1}}.\vadjust{\goodbreak}
\]
For the upper bound, let $H:=-\sum_{n=1}^{\infty}\nu_{n}\log\nu
_{n}$ and let
\[
h :=H/\nu_{\infty}+\log\nu_{\infty}=-\sum_{n=1}^{\infty}\frac
{\nu_{n}}%
{\nu_{\infty}}\log\frac{\nu_{n}}{\nu_{\infty}}.
\]
We need to prove that
\[
\mathcal{P} (\tau>t) \leq\nu_{\infty}e^{(H-t)/\nu_{\infty}} .
\]
If $t$ is such that $\nu_{\infty}e^{(H-t)/\nu_{\infty}}\geq1$, we
are done.
Otherwise, define the sequence of numbers $(\theta_{n})_{n\geq1}$ in
such a way
that for all $n$,
\[
e^{-t\theta_{n}/\nu_{n}}=\nu_{n}e^{(H-t)/\nu_{\infty}}.
\]
The numbers $\theta_{n}$ are positive since $\nu_{n}\leq\nu_{\infty}$,
moreover,
\[
\sum_{n=1}^{\infty}\theta_{n} =\sum_{n=1}^{\infty} \biggl[-\frac1t\nu
_{n}\log
\nu_{n} -\frac1t(H-t)\frac{\nu_{n}}{\nu_{\infty}} \biggr]
=1.
\]
Now,
\begin{eqnarray}
\mathcal{P} (\tau>t) &\leq&\mathcal{P} \Biggl(\bigcup_{n=1}^{\infty}%
\{T_{n}>\theta_{n}t\} \Biggr) \leq\sum_{n=1}^{\infty}\mathcal{P} (T_{n}%
>\theta_{n}t)\nonumber\\[-10pt]\\[-10pt]
&=&\sum_{n=1}^{\infty}e^{-t\theta_{n}/\nu_{n}} =e^{(H-t)/\nu
_{\infty}}\sum
_{n=1}^{\infty}\nu_{n}\nonumber
\end{eqnarray}
and the proof is complete.
\end{pf}

\begin{remark}
Our argument did not make use of the joint law of the $T_{n}$'s which is
unknown. If the r.v.'s $T_{n}$ were independent, a standard exponential bound
would yield $\mathcal{P} (\tau>t) \leq e^{-t/\nu_{1}+h^{\prime}}$,
so we have
the strong feeling that $1/\nu_{\infty}$ is not a sharp bound for the true
rate, which could actually be $1/\nu_{1}$.

It should also be noted that one cannot get rid of $h$ in the previous
statement; one can prove that $\frac{d}{dt}\log\mathcal{P}
(\tau>t) |_{t=0}=0$.
\end{remark}

The following lemma clarifies that either $\tau=\infty$ a.s.\ or
$\tau$
belongs to any interval $[a,b]$ with positive probability.

\begin{lemma}
\label{lemma:tau_support_all_positive}
$\!\!\!$Suppose $\mathcal{P} (\tau
>T)<1$ for
some $T$. Then \mbox{$\mathcal{P} (\tau>a)>\mathcal{P} (\tau>b)$} whenever
$a<b$ and
in particular $\mathcal{P} (\tau>t)<1$ for all $t>0$.
\end{lemma}

\begin{pf}
For $i\geq1$ and $t\geq0$, let
\[
q_{i}(t): =\mathcal{P} (\tau>t|\xi_{0}=i)
=\sum_{n=1}^{\infty}p_{i,n}(t).
\]
By Chapman--Kolmogorov,
%
\begin{equation}
\label{eq:calP_decr}\mathcal{P} (\tau>t) =\sum_{m=1}^{\infty}p_{m}%
(t-s)q_{m}(s) \leq\mathcal{P} (\tau>t-s).\vadjust{\goodbreak}
\end{equation}
Suppose by contradiction that $\mathcal{P} (\tau>a)=\mathcal{P}
(\tau>b)$,
then letting $s=b-a$, $t=b$, we have equality in~\eqref
{eq:calP_decr}. This
implies that $q_{m}(s)=1$ for all $m$ so that we get equality for \emph{any}
choice of $t>0$, meaning that $\mathcal{P} (\tau>\cdot)$ is periodic
as well
as nonincreasing. The only possibility is $\mathcal{P} (\tau>t)=1$
for all $t$.
\end{pf}

We will need the following lemma, that is,\ the probability of no explosion
is larger if one starts from 1 than in any other case. This is a well-known
fact whose proof is a standard exercise we do not repeat here.

\begin{lemma}
\label{lemma_less_expl_from_1} $\mathcal{P} (\tau>t)\leq\mathcal{P}
(\tau>t|\xi_{0}=1)$.
\end{lemma}

\section{Bounds on the energy}

We are finally able to give statements on the decay of the energy
$\mathcal{E}
(t)$ as $t\rightarrow\infty$. First of all we obtain exponential estimates
under $Q$, both in $\mathcal{L}^{1}$ and pathwise. Then we introduce a
smallness condition on $\mathcal{E} (0)$ which is what we need to translate
these results under $P$. The section concludes with the proof of
Theorem~\ref{main_theorem}.

\begin{proposition}
\label{cor_exp_q_ener_expon} Suppose $\nu_{\infty}:=\sum
_{n=1}^{\infty}%
nk_{n}^{-2}<\infty$. Let $X$ be an energy controlled solution and
denote by
$\mathcal{E} (t):=\frac12\sum_{n=1}^{\infty}X_{n}^{2}(t)$ its total
energy at
time $t$. Then
\[
\lim_{t\rightarrow\infty} \mathbb{E} ^{Q}[\mathcal{E} (t)]=0.
\]
In particular there exists a positive number $h$ such that for all
$t\geq0$
%
\begin{equation}
\label{inequality_3}\mathbb{E} ^{Q}[\mathcal{E} (t)]\leq e^{-t/\nu
_{\infty}%
+h}\mathcal{E} (0).
\end{equation}

\end{proposition}

\begin{pf}
By definition
\[
\mathbb{E} ^{Q}[\mathcal{E} (t)] =\frac12\sum_{n} \mathbb{E}
^{Q}[X_{n}%
^{2}(t)] =\mathcal{E} (0)\sum_{n} p_{n}(t),
\]
whence, by Proposition~\ref{prop:uniq_fwd},
\[
=\mathcal{E} (0)\sum_{n} \mathcal{P} (\xi_{t}=n) =\mathcal{E}
(0)\mathcal{P}
(\xi_{t}<\infty) =\mathcal{E} (0)\mathcal{P} (\tau>t).
\]
Finally, we apply Lemma~\ref{lemma_less_expl_from_1} and
Proposition~\ref{prop_tau_exp_bd} and we find~\eqref{inequality_3}.
\end{pf}

Using Borel--Cantelli arguments, one can deduce from (\ref
{inequality_3}) some
$Q$-a.s. statements about the decrease to zero of the energy, at least on
given sequences of times going to infinity. To extend this to all
sequences we
will need the following lemma.

\begin{lemma}
\label{lemma_inequality} Let $X$ be an energy controlled solution and denote
by $\mathcal{E} (t):=\frac12\sum_{n=1}^{\infty}X_{n}^{2}(t)$ its
total energy
at time $t$. Then for every $t\geq s\geq0$ we have%
\[
Q \bigl( \mathcal{E} ( t ) \leq\mathcal{E} ( s )
\bigr) =1.
\]

\end{lemma}

\begin{pf}
Let $s\geq0$ be given and set $\chi=X(s) $. Consider the linear
equation on
$[s,\infty)$ with initial condition $\chi$. It is proved in \cite
{BFM2} that it
has a unique strong solution $Y$, with the property
\[
Q \Biggl( \sum_{n=1}^{\infty}Y_{n}^{2} ( t ) \leq\sum_{n=1}^{\infty
}\chi_{n}^{2} \Biggr) =1
\]
for every $t\in\lbrack s,\infty)$ (the result in \cite{BFM2} is for constant
initial conditions, but the extension to nonanticipative square integrable
random initial conditions is straightforward). But also $X$ restricted to
$[s,\infty)$ is a solution of the same equation, hence, equal to $Y$, on
$[s,\infty)$. The previous identity is thus equal to the claim of the lemma.
\end{pf}

Now we can prove $Q$-a.s.\ exponential decay of energy.

\begin{proposition}\label{prop:Q-as_exp_decay}
Under the same hypothesis of
Proposition \ref{cor_exp_q_ener_expon}, the total energy of solutions
goes to
zero at least exponentially fast pathwise under $Q$,
\[
\limsup_{t\rightarrow\infty}\frac1t\log\mathcal{E} (t)\leq-\frac{1}{\nu_{\infty}},\qquad Q\mbox{-a.s.}
\]

\end{proposition}

\begin{pf}
Let $\epsilon>0$ be given. Set $\alpha:=1/\nu_{\infty}+\epsilon$.
We have
\[
Q \bigl(n^{-1}\log\mathcal{E} (n)> \alpha\bigr)
\leq e^{-\alpha n}\mathbb{E} ^{Q}[\mathcal{E} (n)] \leq Ce^{-\epsilon
n},
\]
where, by Proposition~\ref{cor_exp_q_ener_expon}, $C=e^{h}\mathcal{E} (0)$
does not depend on $n$. Hence, the above probabilities are summable on
$n$ and
by Borel--Cantelli lemma there exists a measurable set $N$ with
$Q(N)=0$ and
the following property: for every $\omega\in N^{c}$ there exists
$n_{0} (
\omega) $ such that, for all $n\geq n_{0} ( \omega) $,
$\mathcal{E} ( n,\omega) \leq e^{-\alpha n}$. Taking the supremum
for $n<n_{0} ( \omega) $, we obtain that there exists a constant
$C ( \omega) >0$ such that
\[
\mathcal{E} ( n,\omega) \leq C ( \omega) e^{-\alpha
n}%
\]
for all $\omega\in N^{c}$ and $n\geq0$. From Lemma \ref{lemma_inequality},
there exists a measurable set $\widetilde{N}$ with $Q(\widetilde
{N})=0$ such
that $\mathcal{E} ( r,\omega) \leq\mathcal{E} ( \lfloor
r\rfloor,\omega) $ for all $\omega\in\widetilde{N}^{c}$ and all
rational numbers $r\in[ 0,\infty) $. This implies, for
$\omega\in N^{c}\cap\widetilde{N}^{c}$,%
\[
\mathcal{E} ( r,\omega) \leq\mathcal{E} ( \lfloor
r\rfloor,\omega) \leq C ( \omega) e^{-\alpha\lfloor
r\rfloor} \leq C^{\prime} ( \omega) e^{-\alpha r }%
\]
with $C^{\prime} ( \omega) =C ( \omega) e^{\alpha
^{\prime}}$, for all $r\in[ 0,\infty) \cap\mathbb{Q}$.

With $Q$ probability one, the function $\mathcal{E} ( t ) $ is
lower semicontinuous, being the supremum in $N$ of the functions $\sum
_{n=1}^{N}X_{n}^{2} ( t ) $ which are continuous. Thus we get%
\[
Q \bigl( \mathcal{E} ( t ) \leq C^{\prime}e^{-\alpha t} \mbox{ for
every }t\geq0 \bigr) =1.
\]
Letting $\epsilon$ go to zero on the rationals completes the\vadjust{\goodbreak} proof.
\end{pf}

The reader should be aware that this proposition does not automatically hold
under $P$. Truly, $P$ and $Q$ are equivalent on all $F_{t}$, but
$C^{\prime}$
is $F_{\infty}$-measurable and not $F_{t}$-measurable for any $t$.

In general, we cannot prove that $P$ and $Q$ are equivalent on
$F_{\infty}$,
so we cannot translate such claim into a similar statement on the original
nonlinear equation (\ref{eq_stratonovich}).

When $\mathcal{E} ( 0 ) $ is small enough, however, we can prove
the equivalence and hence, compute exponential upper bounds for
$\mathcal{E}
(t)$, both $P$-a.s.\ and in mean value.

\begin{proposition}
\label{prop:novikov_infty_time} Let $X$ be an energy controlled
solution and
denote by $\mathcal{E} (t):=\frac12\sum_{n=1}^{\infty}X_{n}^{2}(t)$
its total
energy at time $t$. Suppose $\nu_{\infty}:=\sum_{n=1}^{\infty}nk_{n}
^{-2}<\infty$ and $\nu_{\infty}\mathcal{E} (0)<1$. Then
\[
\mathbb{E} ^{Q} \bigl[e^{\int_{0}^{\infty}\mathcal{E} (t)\,dt} \bigr]<\infty,
\]
so that in particular, $P$ and $Q$ are equivalent on $F_{\infty}$.
\end{proposition}

\begin{pf}
By Proposition~\ref{cor_exp_q_ener_expon} and the definition of energy
controlled solution, we have
\begin{eqnarray*}
\mathbb{E} ^{Q}[\mathcal{E} (t)]&<&\mathcal{E} (0)e^{-t/\nu_{\infty
}+h}\qquad  \forall t\geq0,
\\
0&\leq&\mathcal{E} (t)\leq\mathcal{E} (0),\qquad Q\mbox{-a.s.}
\end{eqnarray*}
Let $X:=\int_{0}^{\infty}\mathcal{E} (t)\,dt\geq0$ and $x\geq0$. Then
\begin{eqnarray*}
xQ(X>x) &\leq&\mathbb{E} ^{Q}[X;X>x] = \int_{0}^{\infty}\mathbb{E}
^{Q}[\mathcal{E} (t);X>x]\,dt\\
&\leq&\int_{0}^{\infty}\min\bigl(\mathbb{E} ^{Q}[\mathcal{E}
(t)];\mathcal{E}
(0)Q(X>x)\bigr)\,dt\\
&\leq&\mathcal{E} (0)\int_{0}^{\infty}\min\bigl(e^{-t/\nu_{\infty
}+h};Q(X>x)\bigr)\,dt\\
&=&\mathcal{E} (0)\int_{0}^{u}Q(X>x)\,dt+\mathcal{E} (0)\int
_{u}^{\infty}%
e^{-t/\nu_{\infty}+h}\,dt,
\end{eqnarray*}
where $u$ is such that $e^{-u/\nu_{\infty}+h}=Q(X>x)$. Hence,
\[
xQ(X>x) \leq\mathcal{E} (0)Q(X>x)u+\mathcal{E} (0)\nu_{\infty
}e^{-u/\nu
_{\infty}+h} =\mathcal{E} (0)Q(X>x)(u+\nu_{\infty}).
\]
If $Q(X>x)=0$ for some $x>0$, then clearly $X$ is bounded and we are done.
Otherwise we get
\[
u\geq\frac x{\mathcal{E} (0)}-\nu_{\infty},
\]
that is
\[
Q(X>x)=e^{-u/\nu_{\infty}+h} \leq e^{- x/{(\nu_{\infty}\mathcal{E}
(0))}+h+1},\vadjust{\goodbreak}
\]
yielding
\[
Q(e^{X}>y) \leq y^{-1/{(\nu_{\infty}\mathcal{E} (0))}}e^{h+1}
\]
and finally
\[
\mathbb{E} ^{Q}[e^{X}]=\int_{0}^{\infty}Q(e^{X}>y)\,dy \leq
1+e^{h+1}\int_{1}^{\infty}y^{-1/{(\nu_{\infty}\mathcal{E} (0))}}\,dy
<\infty,
\]
where we used $\nu_{\infty}\mathcal{E} (0)<1$.
\end{pf}

The heuristic behind the proof, which is quite hidden, that is, given that
$\mathbb{E} ^{Q}[X]$ is bounded, $\mathbb{E} ^{Q}[e^{X}]$ is maximum
if $X$ is
spread as much as possible. Since $X=\int\mathcal{E} (t)\,dt$ and
$\mathcal{E}
(t)\in[0,\mathcal{E} (0)]$, this is done by choosing $\mathcal{E}
(t,\omega)\in\{0,\mathcal{E} (0)\}$, and in particular $\mathcal{E}
(t,\omega)=\mathcal{E} (0)\mathbb{I}_{[0,y(\omega)]}(t)$.

\begin{corollary}
\label{cor:Pas_exp_rate_ener} Under the same hypothesis of
Proposition \ref{prop:novikov_infty_time}, meaning in particular that
$\mathcal{E} (0)<1/\nu_{\infty}$, the energy goes to zero at least
exponentially fast pathwise under $P$,
\[
\limsup_{t\rightarrow\infty}\frac1t\log\mathcal{E} (t) \leq
-\frac1{\nu
_{\infty}}, \qquad P\mbox{-a.s.}
\]
\end{corollary}

\begin{pf}
Just a direct consequence of Propositions~\ref{prop:Q-as_exp_decay}
and \ref{prop:novikov_infty_time}.
\end{pf}

The same condition on the smallness of $\mathcal{E} (0)$ arises when
we want
to establish an exponential decay for the mean value of $\mathcal{E} (t)$
under $P$.

\begin{proposition}
\label{prop:P_avg_exp_rate_ener} If
\[
\mathbb{E} ^{Q}[\mathcal{E} (t)]\leq\mathcal{E} (0)e^{-\alpha t+h}%
\]
then
%
\begin{equation}
\label{eq:exp_bd_P_avg_energy}\mathbb{E} ^{P}[\mathcal{E} (t)]\leq
\mathcal{E}
(0)\exp\bigl( ( 1-1/p )\bigl [h+ \bigl( p\mathcal{E} (0)-\alpha\bigr)
t \bigr] \bigr)
\end{equation}
for every $p>1$. In particular, under the same hypothesis of
Proposition \ref{prop:novikov_infty_time},
%
\begin{equation}
\label{eq:exp_rate_P_avg_energy}\limsup_{t\rightarrow\infty}\frac
1t\log\mathbb{E} ^{P}[\mathcal{E} (t)] \leq-\frac1{\nu_{\infty}} \bigl(
1-\sqrt{\mathcal{E} (0)\nu_{\infty}} \bigr) ^{2}.%
\end{equation}
\end{proposition}

\begin{pf}
The density $f_{t}$ of $P$ with respect to $Q$ on $F_{t}$ is
%
\begin{equation}
\label{density_PQ_time_t}f_{t}=\exp\bigl( M_{t}+\tfrac{1}{2} [ M ]
_{t} \bigr),
\end{equation}
where
\[
M_{t}:=\sum_{n=1}^{\infty}\int_{0}^{t}X_{n} ( s ) \,dW_{n} (
s ) ,\qquad  [ M ] _{t}=\int_{0}^{t}\sum_{n=1}^{\infty}%
X_{n}^{2} ( s ) \,ds.
\]
From (\ref{energy_inequality}) we have
%
\begin{equation}
\exp( \lambda[ M ] _{t} ) \leq e^{2\lambda\mathcal{E}
(0)t} \label{est_1}%
\end{equation}
for every $\lambda>0$.

For every $p,p^{\prime}>1$ with $\frac{1}{p}+\frac{1}{p^{\prime
}}=1$, from the
a.s.\ condition $\mathcal{E} (t)\leq\mathcal{E} (0)$ and the
assumption of the
proposition we have
\begin{eqnarray*}
\mathbb{E} ^{P}[\mathcal{E} (t)] & =&\mathbb{E} ^{Q}[f_{t}\mathcal{E}
(t)]\leq\mathbb{E} ^{Q}[f_{t}^{p}]^{1/p}\mathbb{E} ^{Q}[\mathcal{E}
(t)^{p^{\prime}}]^{1/p^{\prime}}\\
& \leq&\mathbb{E} ^{Q}[f_{t}^{p}]^{1/p}\mathbb{E} ^{Q}[\mathcal{E}
(t)\mathcal{E} (0)^{p^{\prime}-1}]^{1/p^{\prime}}\\
& \leq&\mathcal{E} (0)^{1-1/p^{\prime}}\mathbb{E} ^{Q}[f_{t}^{p}%
]^{1/p}\mathcal{E} (0)^{1/p^{\prime}}e^{-({\alpha}/{p^{\prime
})}t+
h/{p^{\prime}}}\\
& =&\mathcal{E} (0)\mathbb{E} ^{Q}[f_{t}^{p}]^{1/p}e^{-({\alpha}/
{p^{\prime})}t+h/{p^{\prime}}}.
\end{eqnarray*}
From (\ref{density_PQ_time_t}) we have
\begin{eqnarray*}
\mathbb{E} ^{Q}[f_{t}^{p}]&=&\mathbb{E} ^{P}[f_{t}^{p-1}]\\
&=&\mathbb{E}
^{P} \biggl[
\exp\biggl( ( p-1 ) M_{t}+\frac{ ( p-1 ) }{2} [
M ] _{t} \biggr) \biggr]
\\
&=&\mathbb{E} ^{P} \biggl[ \exp\biggl( ( p-1 ) M_{t}-\frac{ (
p-1 ) ^{2}}{2} [ M ] _{t} \biggr) \exp\biggl( \frac{ (
p-1 ) p}{2} [ M ] _{t} \biggr) \biggr]
\end{eqnarray*}
and now we use (\ref{est_1}) to get
\[
\leq e^{ ( p-1 ) p\mathcal{E} (0)t} \mathbb{E} ^{P} \bigl[
e^{ ( p-1 ) M_{t}-({ ( p-1 ) ^{2}}/{2}) [
M ] _{t}} \bigr] =e^{ ( p-1 ) p\mathcal{E} (0)t}
=e^{({p^{2}}/{p^{\prime}})\mathcal{E} (0)t}%
\]
by Girsanov's theorem. To summarize:
\[
\mathbb{E} ^{P}[\mathcal{E} (t)]\leq\mathcal{E} (0)e^{
({1}/{p^{\prime}})p\mathcal{E} (0)t}e^{-({\alpha}/{p^{\prime}})t+ h/{p^{\prime
}}}%
\]
which implies the first claim of the proposition.

To prove the last statement, let $\alpha=1/\nu_{\infty}$. Then
optimization on
$p$ under the condition $\mathcal{E} (0)\nu_{\infty}<1$ gives that
the right-hand side of \eqref{eq:exp_bd_P_avg_energy} is minimum when
$p$ is equal to
$\phi(t)=\sqrt{\frac{1/\nu_{\infty}-h/t}{\mathcal{E} (0)}}$. [We
notice that
$\phi(t)>1$ for $t>\frac h{1/\nu_{\infty}-\mathcal{E} (0)}$.] With
a simple
computation, \eqref{eq:exp_bd_P_avg_energy} becomes
\[
\mathbb{E} ^{P}[\mathcal{E} (t)]\leq\mathcal{E} (0)\exp\biggl\{ -(p-1) \biggl(
\frac{\phi^{2}(t)}p-1 \biggr) \mathcal{E} (0)t \biggr\}.
\]
Letting $p=\phi(t)>1$, we get
\[
\frac1t\log\mathbb{E} ^{P}[\mathcal{E} (t)] \leq\frac1t\log
\mathcal{E}
(0)-\bigl(\phi(t)-1\bigr)^{2}\mathcal{E} (0),\qquad  t>\frac h{1/\nu_{\infty
}-\mathcal{E}
(0)}.
\]
Taking the limsup for $t\rightarrow\infty$ leads to \eqref
{eq:exp_rate_P_avg_energy}.
\end{pf}

\subsection{\texorpdfstring{Proof of Theorem \protect\ref{main_theorem}}{Proof of Theorem 1}}
\label{proof main theorem}

We have $\nu_{\infty}=\sum_{n=1}^{\infty}nk_{n}^{-2}<\infty$, so
by virtue of
Proposition \ref{prop:A_not_regular}, $\mathcal{P} (\tau>t)<1$ for
some $t$.
Then by Lemma \ref{lemma:tau_support_all_positive}, $\mathcal{P}
(\tau>t)<1$
for all $t$. By equation \eqref{eq:Qavg_ener_via_tau} $\mathbb{E}
^{Q}[\mathcal{E} (t)]<\mathcal{E} (0)$ and since $\mathcal{E}
(t)\leq
\mathcal{E} (0)$ $Q$-a.s.\ we get $Q(\mathcal{E} (t)=\mathcal{E}
(0))<1$, for
all $t$. Equivalence of $P$ and $Q$ on $F_{t}$ yields the first statement.

By Proposition \ref{cor_exp_q_ener_expon}, $Q(\mathcal{E}
(t)>\epsilon)\leq
Ke^{-t/\nu_{\infty}}$ for some $K$ not depending on~$t$, hence, for
$t$ large
enough, this event has $Q$- (and hence $P$-) probability less than~1,
so we
proved the second statement.

The third statement is proved in Corollary \ref{cor:Pas_exp_rate_ener} and
Proposition \ref{prop:P_avg_exp_rate_ener}.

The proof of Theorem \ref{main_theorem} is complete.

\section{Lack of regular solutions}\label{sec:lack_reg}

As a consequence of our result on the dissipation of energy we can
prove that there exists no regular solution.

Define the space
\[
V= \Biggl\{ x\in l^{2}\dvtx\sum_{n=1}^{\infty}k_{n}^{2}x_{n}%
^{2}<\infty\Biggr\}
\]
which is an Hilbert space under the norm $ \Vert x \Vert
_{V}^{2}=\sum_{n=1}^{\infty}k_{n}^{2}x_{n}^{2}$.

\begin{proposition}
Assume that $ \{ X ( t ) ;t\in[ 0,T ] \} $
is an energy controlled solution. Then,
\[
P \biggl( \int_{0}^{T} \Vert X ( t ) \Vert_{V}%
^{2}\,dt=\infty\biggr) >0.
\]
\end{proposition}

\begin{pf}
We will actually prove the following statement. Assume that $ \{
X ( t ) ; t\in[ 0,T ] \} $ is an energy
controlled solution such that
%
\begin{equation}\label{eq:regularity_condition}
P \biggl( \int_{0}^{T} \Vert X ( t ) \Vert_{V}%
^{2}\,dt<\infty\biggr) =1.
\end{equation}
Then, for every $t\in[ 0,T ] $, $P( \energy(t) =\energy(0))
=1$.

Since the latter is in contradiction with Theorem \ref{main_theorem},
then \eqref{eq:regularity_condition} will be proven to be false.

Let $X$ be a solution as in the claim. By It\^{o}'s formula,
\begin{eqnarray*}
d \Biggl( \sum_{n=1}^{N}X_{n}^{2} \Biggr) & =&2\sum_{n=1}^{N} (
k_{n-1}X_{n-1}^{2}X_{n}-k_{n}X_{n}^{2}X_{n+1} ) \,dt\\
&&{} -\sum_{n=1}^{N} ( k_{n}^{2}+k_{n-1}^{2} ) X_{n}^{2}\,dt
+\sum_{n=1}^{N} ( k_{n-1}^{2}X_{n-1}^{2}+k_{n}^{2}X_{n+1}^{2} ) \,dt\\
&&{} +2\sum_{n=1}^{N} ( k_{n-1}X_{n-1}X_{n}\,dW_{n-1}-k_{n}X_{n}%
X_{n+1}\,dW_{n} )\\
& =&-2k_{N}X_{N}^{2}X_{N+1}\,dt\\
&&{} -k_{N}^{2}X_{N}^{2}\,dt+k_{N}^{2}X_{N+1}^{2}\,dt\\
&&{} -2 k_{N}X_{N}X_{N+1}\,dW_{N}.
\end{eqnarray*}
Hence,
\begin{eqnarray}\label{eq:resti}
\sum_{n=1}^{N}X_{n}^{2} ( t ) -\sum_{n=1}^{N}x_{n}^{2}
&=&\int_{0}^{t} (-2k_{N}X_{N}^{2}X_{N+1}+ k_{N}^{2}%
X_{N+1}^{2}-k_{N}^{2}X_{N}^{2} ) \,ds\nonumber\\[-8pt]\\[-8pt]
&&{}-\int_{0}^{t}2 k_{N}X_{N}X_{N+1}\,dW_{N} ( s ) .\nonumber
\end{eqnarray}
Notice that $|X_{N+1}|\leq\Vert x\Vert_{l^2} $, $k_{N}\geq1$ and
$k_{N}\leq k_{N+1}$, so that the first integral can be bounded by
\[
C\int_{0}^{t} (k_{N}^2X_{N}^{2}+ k_{N+1}^{2}
X_{N+1}^{2}+k_{N}^{2}X_{N}^{2} )\, ds.
\]
Then the a.s.\ inequality $\int_{0}^{T} \Vert X ( t )
\Vert_{V} ^{2}\,dt<\infty$ implies that
$\int_{0}^{T}k_{N}^{2}X_{N}^{2}\,ds$ goes to zero a.s.\ and in
probability.

The latter is true also for the stochastic integral in
\eqref{eq:resti} by continuity in probability of stochastic integrals;
since $k_NX_NX_{N+1}$ converges to zero in probability in $L^2(0,T)$,
its stochastic integral converges to zero in probability. This can be
easily proved by applying the following well-known inequality of
stochastic integrals. For every $\gamma,\delta>0$,
\[
P \biggl( \biggl\vert\int_{0}^{t}k_{N}X_{N}X_{N+1}\,dW_{N} ( s )
\biggr\vert>\gamma\biggr) \leq P \biggl( \int_{0}^{t}k_{N}^{2}X_{N}%
^{2}X_{N+1}^{2}\,ds>\delta\biggr) +\frac{\delta}{\gamma^{2}}.
\]

Finally, all terms on the RHS of \eqref{eq:resti} converge to zero in
probability. We get
\[
\sum_{n=1}^{N}X_{n}^{2} ( t )
\mathop{\longrightarrow}_{N\rightarrow\infty}^{P} \Vert x \Vert_{l^{2}}^{2}.
\]
But we know, by $P$-a.s. monotonicity, that
\[
\sum_{n=1}^{N}X_{n}^{2} ( t )
\mathop{\longrightarrow}_{N\rightarrow\infty}^{P}\sum_{n=1}^{\infty}X_{n}^{2}
( t ) .
\]
Hence, the claim is proved and the proposition as well.
\end{pf}

\begin{remark}
If one could prove local existence (up to some random time) of
solutions in $V$ with initial conditions in $V$, the above proposition
would show a blow-up in the $V$ norm. This appears to be a difficult
open problem in the stochastic case, while in the deterministic case,
it is well known (see, e.g., \cite{Ces}).
\end{remark}


%

\printaddresses

\end{document}